\documentclass[12pt]{article}

\usepackage{amsmath}
\usepackage{amssymb}
\usepackage{amscd}
\usepackage[all]{xy}

\newtheorem{theorem}{Theorem}[section]

\newtheorem{proposition}[theorem]{Proposition}
\newtheorem{corollary}[theorem]{Corollary}
\newtheorem{lemma}[theorem]{Lemma}

\newenvironment{proof}{\begin{trivlist} \item[]{\bf Proof.}}
{\par\hfill $\square$\end{trivlist}}

\newcommand{\Lone}{{{\rm L}^1}}

\newcommand{\Linfty}{{{\rm L}^\infty}}
\newcommand{\C}{\mathbb{C}}
\newcommand{\N}{\mathbb{N}}
\renewcommand{\P}{\mathbb{P}}
\newcommand{\R}{\mathbb{R}}

\newcommand{\K}{{\cal K}}

\newcommand{\ddc}{{{\rm dd}^{\rm c}}}
\renewcommand{\d}{{\rm d}}

\newcommand{\supp}{{\rm supp}}

\title{Decay of correlations for H\'enon maps}
\author{Tien-Cuong Dinh}

\begin{document}

\maketitle

\begin{abstract} We show, for a class of automorphisms of $\C^k$, that their 
equilibrium measures are exponentially mixing. In particular, this holds for (generalized) H\'enon maps in $\C^2$.
\end{abstract}

\section{Introduction}

Let $f$ be a polynomial automorphism of $\C^k$. We also write $f$ for its extension as a birational self-map
of $\P^k$. When the exceptional sets $I_\pm$ of $f^\pm$ are not empty and 
satisfy $I_+\cap I_-=\emptyset$,
we say that $f$ is a {\it regular automorphism} in the sense of Sibony. The reader may find a description 
of these maps in the survey of Sibony \cite{Sibony}. See also Friedland-Milnor \cite{FriedlandMilnor}, 
Bedford-Lyubich-Smillie \cite{BedfordSmillie, BedfordLyubichSmillie} and Forn\ae ss-Sibony \cite{FornaessSibony1}.
We recall here some facts. 

The exceptional sets $I_\pm$ are contained in the 
hyperplane at infinity $L_\infty$. There
exists an integer $s$ such that $\dim I_+=k-1-s$ and $\dim I_-=s-1$. 
We have  $f(L_\infty\setminus I_+)=I_-$ and $f^{-1}(L_\infty\setminus I_-)=I_+$. 
Moreover, $I_-$ is attractive for $f$ and $I_+$ is attractive for $f^{-1}$. 
If $d_+$ and $d_-$ are the
algebraic degrees of $f$ and $f^{-1}$ respectively, then $d_+^s=d_-^{k-s}>1$. 
We have $d_+=d_-$ when $k=2s$.
In the dimension $k=2$, 
$f$ is a generalized H\'enon automorphism which are the only dynamically 
interesting polynomial automorphisms of $\C^2$.

Sibony constructed for such a map an invariant 
probability measure $\mu$ (called {\it Green measure} or
{\it equilibrium measure}) as 
an exterior product of positive closed $(1,1)$-currents:
$$\mu= T_+^s \wedge T_-^{k-s}$$
where $T_\pm$ are {\it Green currents of bidegree $(1,1)$} and of mass 1 associated to $f^{\pm 1}$.
They have local continuous potentials in $\P^k\setminus I_\pm$ and satisfy $f^*(T_+)=d_+T_+$,
$f_*(T_-)=d_-T_-$. The current $T_+^s$
(resp. $T_-^{k-s}$)  is supported in the boundary of the 
{\it filled Julia set} $\K_+$ (resp. $\K_-$); it is called 
{\it Green current of bidegree $(s,s)$} (resp. $(k-s,k-s)$)
associated to $f$ (resp. to $f^{-1}$).

Recall that $\K_+$ (resp. $\K_-$) is the set of points $z\in\C^k$ such that
the orbit $(f^n(z))_{n\in\N}$ (resp. $(f^{-n}(z))_{n\in\N}$) is bounded in $\C^k$. We have
$\overline\K_\pm\cap L_\infty=I_\pm$. 
The open set $\P^k\setminus\overline \K_+$ (resp. $\P^k\setminus\overline \K_-$) is the immediate bassin 
of $I_-$ for $f$ (resp. $I_+$ for $f^{-1}$).
The measure $\mu$ is supported in the boundary of the compact set $\K:=\K_+\cap \K_-$. 

It was recently proved in \cite{DinhSibony5, Guedj} 
that $\mu$ is mixing. 
This generalizes results of Bedford-Smillie \cite{BedfordSmillie} and Sibony \cite{Sibony}.
The proofs follow the same approach and use the property that $T_+^s$ and 
$T_-^{k-s}$ are extremal currents. 
In this paper, we use another method to show that $\mu$ is mixing and the speed of mixing is exponential
when $k=2s$. We keep the above notation.

\begin{theorem}
Let $f$ be as above and assume that $k=2s$. 
Then, there exists a constant $c>0$ such that 
$$\left|\int (\varphi\circ f^n)\psi\d\mu -\left(\int \varphi\d\mu\right) \left(\int \psi\d\mu\right)\right|
\leq cd_+^{-n/2}\|\varphi\|_{{\cal C}^2}\|\psi\|_{{\cal C}^2}$$
for all $n\geq 0$ 
and all real-valued ${\cal C}^2$ functions $\varphi$ and $\psi$ in $\C^k$.
\end{theorem} 

Of course, this result holds for polynomial automorphisms of positive entropy in $\C^2$, 
in particular, for H\'enon maps. We can apply it for any real H\'enon map of degree 
$d$ in $\R^2$ 
which admits an invariant probability measure $\mu$ of entropy $\log d$.

In \cite{BedfordLyubichSmillie}, Bedford-Lyubich-Smillie proved for complex 
H\'enon maps that the equilibrium measure
is Bernoulli. This is the strongest mixing in the sense of measures. However, it does not imply the
decay of correlations in our sense. 

Observe also that in Theorem 1.1 we cannot replace $\|\varphi\|_{{\cal C}^2}$
or $\|\psi\|_{{\cal C}^2}$ by $\|\varphi\|_{\Linfty}$ or $\|\psi\|_{\Linfty}$ since, in general, 
$\psi\circ f^{-n}$ and $\varphi\circ f^n$ do not converge in $\Lone(\mu)$ to a constant.

We think that Theorem 1.1 should be true for every regular automorphism, 
for H\'enon-like maps in dimension 2 \cite{Dujardin,
DinhDujardinSibony} and for 
larger classes of birational maps of $\P^k$ considered in \cite{DinhSibony5}.
For a general regular automorphism, 
in order to apply our approach, one needs to analyse the indeterminacy sets of some
automorphisms close to the regular polynomial maps in the sense of \cite{DinhSibony2} (see Lemma 3.2 below). 

Note that the exponential decay of correlations has been proved for some polynomial-like maps and for 
meromorphic maps of large topological degree in \cite{FornaessSibony2, DinhSibony1, DinhSibony3}. 
In these cases, by the classical Gordin-Liverani theorem, our estimates imply the central limit theorem
for bounded quasi-p.s.h. observables.

In Sections 2 and 3, we give some properties of the Green currents and the equilibrium measure. 
The method of $\ddc$-resolution developed in \cite{DinhSibony1, DinhSibony3,
DinhSibony4, DinhSibony5} will be applied to establish the 
necessary estimates (Propositions 2.1, 3.1). 
We then deduce in Section 4 
the mixing and the speed of mixing. 

\section{Convergence toward the Green current}

Let us recall two properties of currents on $\P^k$ that will be used later on.
Since $\P^k$ is homogeneous, every positive closed 
current $S$ on $\P^k$ can be regularized on every neighbourhood $U$ 
of $\supp(S)$. If $T$ is a positive closed $(1,1)$-current with local continuous potentials in a neighbourhood
of $\overline U$, then the positive closed current $T^m\wedge S$ is well defined and depends continuously
on $S$. We refer to  \cite{Federer, Lelong, BedfordTaylor, Demailly, Sibony} for the basics of
the theory of currents.

Now, consider a regular automorphism $f$ on $\C^k$ as in Section 1. We do not suppose that $k=2s$. 
Fix neighbourhoods 
$U_i$ of $\overline \K_+$ and $V_i$ of $\overline \K_-$ 
such that $f^{-1}(U_i)\Subset U_i$, $U_1\Subset U_2$, $f(V_i)\Subset V_i$, $V_1\Subset V_2$
and $U_2\cap V_2\Subset \C^k$. Observe that $\K_+\cap \K_-\subset U_1\cap V_1$. 

Let $\Omega$ be a real $(k-s+1,k-s+1)$-current with support in $\overline V_1$. 
Assume that there exists a positive closed $(k-s+1,k-s+1)$-current $\Omega'$ supported in $\overline V_1$
such that $-\Omega'\leq \Omega\leq \Omega'$.
Define the norm $\|\Omega\|_*$ of $\Omega$ as
$$\|\Omega\|_*:=\min\{\|\Omega'\|,\ \Omega' \mbox{ as above} \}$$
where $\|\Omega'\|=\langle\Omega',\omega^{s-1}\rangle$ is the mass of $\Omega'$.
Here $\omega$ denotes the Fubini-Study form on $\P^k$ normalized by $\int \omega^k=1$.
\\

Keeping the above notation, the main result of this section is the following proposition.

\begin{proposition}
Let $R$ be a positive closed $(s,s)$-current 
of mass $1$ supported in $U_1$ and smooth on $\C^k$. Let $\Phi$ be 
a real smooth $(k-s,k-s)$-form with compact support in $V_1\cap \C^k$. Assume that 
$\ddc \Phi\geq 0$ in $U_2$. 
Then, there exist constants $c>0$ independent of $R$, $\Phi$,  and $c_R>0$ independent of $\Phi$ such that
$$\langle d_+^{-sn} f^{n*}(R) - T_+^s,\Phi\rangle \leq cd_+^{-n} \|\ddc\Phi\|_*$$
and
$$\left|\langle d_+^{-sn} f^{n*}(R) - T_+^s,\Phi\rangle \right|\leq c_Rd_+^{-n} \|\ddc\Phi\|_*$$
for every $n\geq 0$. 
In particular, $d_+^{-sn}f^{n*}(R)\rightarrow T_+^s$ as $n\rightarrow\infty$.
\end{proposition} 

The current $f^{n*}(R)$ is well defined since $f^{-n}$ is holomorphic in $U_1$.
We have
\begin{eqnarray}
\langle d_+^{-sn} f^{n*}(R) -T_+^s,\Phi\rangle & = & d_+^{-sn} \langle f^{n*}(R-T_+^s), \Phi \rangle \nonumber\\
& = & d_+^{-sn} \langle R-T_+^s, (f^n)_*\Phi \rangle.
\end{eqnarray}
Since the currents $R$ and $T_+^s$ have the same mass 1, they are cohomologous. 
On $\P^k$, 
$R-T_+^s$ is $\ddc$-exact.
Hence, the last term in (1)
does not change if we subtract a $\ddc$-closed form from $(f^n)_*\Phi$.
We will use the following
lemma applied to $\ddc(f^n)_*\Phi$.

\begin{lemma} Let $\Omega$ be a real smooth form of bidegree
$(k-s+1,k-s+1)$ supported in $\overline V_1$ such that 
$\Omega\geq 0$ on $U_2$ and $\|\Omega\|_*\leq 1$. Assume that $\Omega$ is $\ddc$-exact. 
Then there exist $c>0$ independent of
$\Omega$ and a real continuous 
$(k-s,k-s)$-form $\Psi$ such that $\ddc \Psi=\Omega$, $\|\Psi\| \leq c$,
$\Psi\leq 0$ on $U_1$ and $\Psi\geq -c\omega^{k-s}$ on $\P^k\setminus V_2$.
\end{lemma}
\begin{proof} 
By Hodge theory \cite{GriffithsHarris}, we have
$$H^{k,k}(\P^k\times \P^k,\C)\simeq \sum_{p+p'=k} H^{p,p}(\P^k,\C) \otimes
H^{p',p'}(\P^k,\C).$$
Hence, 
if $\Delta$ is the diagonal of $\P^k\times \P^k$,
there exists a smooth real $(k,k)$-form $\alpha(x,y)$ on $\P^k\times \P^k$,
cohomologous to $[\Delta]$, with
$\d_x\alpha=\d_y\alpha=0$.
Since $\P^k\times\P^k$ is homogeneous, following 
\cite[Prop. 6.2.3]{BostGilletSoule} (see also \cite{GilletSoule, DinhSibony4, DinhSibony5}), 
one can construct a negative $(k-1,k-1)$-form $K(x,y)$ on 
$\P^k\times \P^k$, smooth outside $\Delta$, such that 
$\ddc K=[\Delta]-\alpha$ and $|K(x,y)|\leq A|x-y|^{1-2k}$ for some constant $A>0$.
Here $|x-y|$ denotes the distance between $x$ and $y$.

Define
$$\Psi'(x):=\int_y K(x,y)\wedge \Omega(y).$$
From the estimate,
one check easily that $\Psi'$ is continuous and $\|\Psi'\|\leq c'$, $\Psi'\leq c'\omega^{k-s}$ on $U_1$, 
$\Psi'\geq -c'\omega^{k-s}$ on $\P^k\setminus V_2$, where $c'>0$ is independent of $\Omega$. 
Define $\Psi:=\Psi'-c'\omega^{k-s}$. 
We obtain  $\|\Psi\|\leq 2c'$, $\Psi\leq 0$ on $U_1$ and 
$\Psi\geq -2c'\omega^{k-s}$ on $\P^k\setminus V_2$. 
We only have to verify that $\ddc\Psi'=\Omega$.

Since $\Omega$ is $\ddc$-exact and $\d_x\alpha=\d_y\alpha=0$, we have 
\begin{eqnarray*}
\ddc\Psi'(x) & := & \int_y (\ddc)_x K(x,y)\wedge \Omega(y)
 =  \int_y \ddc K(x,y)\wedge\Omega(y)\\
& = & \int_y ([\Delta]-\alpha)\wedge \Omega(y)
 =  \Omega(x) -\int_y\alpha\wedge \Omega(y)\\
& = & \Omega(x).
\end{eqnarray*}
Hence, $\ddc\Psi=\ddc\Psi'=\Omega$.
\end{proof}

\noindent
{\bf Proof of Proposition 2.1.} We can assume that $\|\ddc\Phi\|_*=1$. The constants $c$ and $c_i$  
below are independent of $\Phi$ and $R$.  
Define $\Omega:=\ddc\Phi$. By hypotheses, there exists a positive closed current $\Omega'$ of mass 1 supported 
in $\overline V_1$ such that
$-\Omega'\leq \Omega\leq \Omega'$.  
Define $\Omega_n:=\ddc (f^n)_*\Phi= (f^n)_*\Omega$ and 
$\Omega_n':=(f^n)_*\Omega'$. 
These currents have supports in $\overline V_1$ since $f^n(V_1)\Subset V_1$. 
We also have 
$-\Omega_n'\leq\Omega_n\leq \Omega_n'$ and $\Omega_n\geq 0$ on $U_2$ since 
$f^{-n}(U_2)\Subset U_2$.  
A simple calculus on cohomology gives 
$\|\Omega_n'\|=d_+^{(s-1)n}\|\Omega'\|=d_+^{(s-1)n}$.
Lemma 2.2 implies the existence of $\Psi_n$ cohomologous to $(f^n)_*\Phi$ such that 
$\Psi_n\leq 0$ on $U_1$, $\Psi_n\geq -cd_+^{(s-1)n}\omega^{k-s}$ on $\P^k\setminus V_2$ 
and $\|\Psi_n\|\leq c d_+^{(s-1)n}$. In particular, $\Psi_n\leq 0$ on $\supp(R)$. Therefore, 
we deduce from (1) that
\begin{eqnarray}
\langle d_+^{-sn}f^{n*}(R)-T_+^s,\Phi\rangle = d_+^{-sn}\langle R-T_+^s,\Psi_n\rangle  \leq 
 -d_+^{-sn}\langle T_+^s,\Psi_n\rangle.
\end{eqnarray}

We have to bound $-\langle T_+^s,\Psi_n\rangle$. 
Since $T_+$ has continuous potentials in $\P^k\setminus I_+$, we can  
write $T_+ =\omega +\ddc u$ with $u\leq 0$ and $u$ continuous on $\P^k\setminus I_+$. 
One has
\begin{eqnarray}
|\langle T_+^s,\Psi_n\rangle| & = & |\langle \omega\wedge T_+^{s-1} +\ddc (uT_+^{s-1}), \Psi_n \rangle| \nonumber\\
& \leq & |\langle T_+^{s-1},\omega\wedge \Psi_n \rangle| + |\langle uT_+^{s-1}, 
\ddc\Psi_n\rangle|\nonumber\\
& \leq & |\langle T_+^{s-1},\omega\wedge\Psi_n\rangle| -\langle uT_+^{s-1}, \Omega_n'\rangle.
\end{eqnarray}
Since $\Omega_n'$ has support in $\overline V_1$ where $u$ is bounded, the second term
in the last line of (3) is dominated by 
$c_1\langle T_+^{s-1},\Omega_n'\rangle$. The integral $\langle T_+^{s-1},\Omega_n'\rangle$ is cohomological; it
is equal to $\|\Omega_n'\|$. Hence, $-\langle u T_+^{s-1},\Omega_n'\rangle \leq c_1 d_+^{(s-1)n}$. 

For the first term in the last line of (3), 
we write $T_+^{s-1}=\omega\wedge T_+^{s-2}+\ddc(uT_+^{s-2})$. Using expansions
as in (3) and an induction argument, we get $|\langle T_+^{s-1},\omega\wedge \Psi_n\rangle| \leq 
c_2d_+^{(s-1)n}$. At the last step of the induction, we use the inequality $\|\Psi_n\|\leq c d_+^{(s-1)n}$. 
Hence, the first part of Proposition 2.1 follows.
\\

For the second part, 
it is sufficient to prove that 
$|\langle R,\Psi_n\rangle|\leq c_R'd_+^{(s-1)n}$ with $c_R'$ independent of $\Phi$. 
This follows directly from the smoothness of $R$ on $\C^k$
and the properties that  $\|\Psi_n\|\leq cd_+^{(s-1)n}$ and
$-cd_+^{(s-1)n}\omega^s\leq\Psi_n\leq 0$ 
on the neighbourhood $U_1\setminus V_2$ of the singularities of $R$.
\\

Now, we show that $d_+^{-sn}f^{n*}(R)\rightarrow T_+^s$ on $\C^k$.
Consider a real smooth 
test $(k-s,k-s)$-form $\Phi$ with compact support in $\C^k$. 
We want to prove that $\langle d^{-sn}f^{n*}(R)-T_+^s,\Phi\rangle \rightarrow 0$.
Observe that $\P^k\setminus I_+$ is a
union of compact algebraic sets of dimension $s$. Hence, we can construct a positive closed 
$(k-s,k-s)$-form $\Theta$ supported in $\P^k\setminus I_+$ and strictly positive on $\supp(\Phi)$. 
Since 
$$\langle d_+^{-sn} f^{n*}(R)-T_+^s,\Phi\rangle = \langle d_+^{-s(n-m)} f^{(n-m)*}(R)-T_+^s, d_+^{-sm}(f^m)_*\Phi
\rangle,$$
replacing $\Phi$ and $\Theta$ by $d_+^{-sm}(f^m)_*\Phi$ and $(f^m)_*\Theta$, $m$ big enough, one can assume that 
$\supp(\Theta)\subset V_1$. 

Consider a smooth function $\chi$ with compact support in $\C^k$ and strictly p.s.h. on neighbourhood of
$U_2\cap V_2$. Write $\Phi=(\Phi+A\chi\Theta)-A\chi\Theta$ with $A>0$ big enough, so that
$\ddc (\Phi+A\chi\Theta)$ and $\ddc(A\chi\Theta)$ are positive on $U_2$. Hence, it is sufficient to
consider the case where $\ddc\Phi\geq 0$ on $U_2$. 
The second part of the proposition implies that $\langle d^{-sn}f^{n*}(R)-T_+^s,\Phi\rangle \rightarrow 0$.

\hfill $\square$

\section{Convergence toward the Green measure}

In this section, we will apply Proposition 2.1 to the automorphism $F$ constructed in Lemma
3.2 below. Our main result is the following proposition.

\begin{proposition} Let $f$ be as above with $k=2s$. 
Let $\varphi$ be a smooth function on $\P^k$ and p.s.h. on $U_2\cap V_2$. 
Let $R$ (resp. $S$) be a positive closed $(s,s)$-current of mass $1$ 
with support in $U_1$ (resp. in $V_1$) and smooth on $\C^k$.
Then, there exist constants $c>0$ independent of $\varphi$, $R$, $S$, and $c_{R,S}>0$
independent of $\varphi$ such that 
$$\left\langle d_+^{-2sn}f^{n*}(R)\wedge (f^n)_*(S) - \mu, \varphi
\right\rangle \leq cd_+^{-n} \|\varphi\|_{{\cal C}^2}$$
and
$$\big|\left\langle d_+^{-2sn}f^{n*}(R)\wedge (f^n)_*(S) - \mu, \varphi
\right\rangle \big| \leq c_{R,S}d_+^{-n} \|\varphi\|_{{\cal C}^2}$$
for every $n\geq 0$. In particular, $d_+^{-2sn}f^{n*}(R)\wedge (f^n)_*(S)\rightarrow\mu$
as $n\rightarrow\infty$.
\end{proposition}

We will use $z$, $w$ and $(z,w)$ for the canonical coordinates of complex spaces $\C^k$ and
$\C^k\times\C^k$. 
Consider also the canonical inclusions of $\C^k$ and $\C^k\times\C^k$ in $\P^k$ and $\P^{2k}$.
We write $[z:t]$, $[w:t]$ or $[z:w:t]$ for the homogeneous coordinates of projective spaces.
The hyperplanes at infinity are defined by $t=0$. 
If $g:\C^k\rightarrow\C^k$ is a polynomial automorphism, 
we write $g_h$ (resp. $g^{-1}_h$) for the homogeneous part of maximal degree of $g$
(resp. of $g^{-1}$). They are  self-maps of $\C^k$, not invertible in general. 
In the sequence, we always assume that  $k=2s$.

\begin{lemma} Let $F$ be the automophism
of $\C^k\times\C^k$ defined by $F(z,w):=(f(z), f^{-1}(w))$. Then $F$ is regular. 
The indeterminacy sets $I^F_\pm$ of $F^\pm$ are defined by
$$I^F_\pm:=\left\{[z:w:0],\ f^{\pm 1}_h(z)=0,\ f^{\mp 1}_h(w)=0 \right\}.$$
Moreover,
if $\Delta:=\{z=w\}$ is the diagonal of $\C^k\times\C^k$, then $I^F_\pm$ 
do not intersect $\overline \Delta$. In particular, $F(\overline \Delta)\cap \{t=0\} \subset I^F_-$.   
\end{lemma}
\begin{proof} Since $k=2s$, we have $d_+=d_-$ and $F^{\pm 1}_h(z,w)=(f^{\pm 1}_h(z),f^{\mp 1}_h(w))$.
It follows that 
$$I^F_\pm=\left \{ [z:w:0],\ F_h^{\pm 1}(z,w)=0\right\}= 
\left\{[z:w:0],\ f^{\pm 1}_h(z)=f^{\mp 1}_h(w)=0 \right\}.$$
We also have
$$I_\pm:=\left\{ [z:0],\ f^{\pm 1}_h(z)=0\right\}$$
and since $f$ is regular,  
$$\left\{z\in\C^k,\  f_h(z)= f^{-1}_h(z)=0\right\}= \{0\}.$$
This implies that $I^F_+\cap I^F_-=\emptyset$. Hence, $F$ is regular. We also have
$$I^F_\pm \cap \overline\Delta = \left\{ [z:z:0],\ f_h(z)=f^{-1}_h(z)=0\right \}=\emptyset.$$
\end{proof}

\begin{lemma} Under the notation of Lemma 3.2, 
the Green current of bidegree $(2s,2s)$ of $F$ is equal to $T_+^s\otimes T_-^s$.
\end{lemma}
\begin{proof} Let $R$ and $S$ be as in Proposition 3.1. 
Replacing $R$ and $S$ by $d_+^{-s}f^*(R)$ and $d_+^{-s}f_*(S)$, we get
$\supp(R)\cap \{t=0\}\subset I_+$ and $\supp(S)\cap \{t=0\}\subset I_-$. 

Consider the current $R\otimes S$ in $\C^k\times \C^k$ and in $\P^{2k}$. Lemma 3.2 implies  
$\overline{\supp(R\otimes S)}\cap \{t=0\}\subset I_+^F$. Since 
$\dim I_+^F=2s-1$, the trivial extension of $R\otimes S$
in $\P^{2k}$ (that we denote also by $R\otimes S$) is a positive closed current \cite{HarveyPolking}. 
One can check that the mass of $R\otimes S$ is equal to 1.  Proposition 2.1 applied to $F$ 
implies that $d_+^{-2sn} F^{n*}(R\otimes S)$
converge to the Green current of bidegree $(2s,2s)$ of $F$. On the other hand, we have
$$d_+^{-2sn}F^{n*}(R\otimes S)=d_+^{-2sn}f^{n*}(R)\otimes (f^n)_*(S)\rightarrow T_+^s\otimes T_-^s$$ 
in $\C^k\times\C^k$.
Hence, $T_+^s\otimes T_-^s$ is the Green $(2s,2s)$-current of $F$.
\end{proof}

\noindent
{\bf Proof of Proposition 3.1.} 
We can assume that $\varphi$ has compact support in $\C^k$ and $\|\varphi\|_{{\cal C}^2}=1$.
As in Lemma 3.3, we can assume that
the current $R\otimes S$
in $\P^{2k}$ satisfies
$\supp(R\otimes S)\cap \{t=0\}\subset I^F_+$. 

Define $\widehat \varphi(z,w):=\varphi(z)$. 
Since $T_\pm$ are invariant and have continuous potentials out of $I_\pm$, we can write
\begin{eqnarray}
\lefteqn{
\left\langle d_+^{-2sn}f^{n*}(R)\wedge (f^n)_*(S) - \mu, \varphi
\right\rangle} \nonumber\\ 
& = & \left\langle d_+^{-2sn} f^{n*}(R)\otimes (f^n)_*(S) - T_+^s\otimes T_-^s, \widehat\varphi[\Delta] 
\right\rangle\nonumber.
\end{eqnarray}
Using a regularization of $[\Delta]$, 
one may find a smooth current $\Theta$ of mass 1 supported in a small neighbourhood
${\cal W}$ of $\overline\Delta$  such that 
\begin{eqnarray}
\lefteqn{\big|\left\langle d_+^{-2sn} f^{n*}(R)\otimes (f^n)_*(S) - T_+^s\otimes T_-^s, \widehat\varphi[\Delta] 
\right\rangle -}\nonumber \\
&& -\left\langle d_+^{-2sn} f^{n*}(R)\otimes (f^n)_*(S) - T_+^s\otimes T_-^s, \widehat\varphi\Theta 
\right\rangle  \big|\leq  d_+^{-n}.\nonumber
\end{eqnarray}
The current $\Theta$ depends on $n$ and ${\cal W}\cap I^F_+=\emptyset$ (see Lemma 3.2).

We have to estimate
$$\left\langle d_+^{-2sn} f^{n*}(R)\otimes (f^n)_*(S) - T_+^s\otimes T_-^s, \widehat\varphi\Theta 
\right\rangle.$$
Fix an integer $m>0$ big enough. Write
\begin{eqnarray*}
\lefteqn{\left\langle d_+^{-2sn} f^{n*}(R)\otimes (f^n)_*(S) - T_+^s\otimes T_-^s, \widehat\varphi\Theta 
\right\rangle } \\
& = & \left\langle d_+^{-2sn} F^{n*}(R\otimes S) - d_+^{-2sm} F^{m*}(T_+^s\otimes T_-^s), 
\widehat\varphi \Theta\right\rangle\\
& = & \left\langle d_+^{-2s(n-m)} F^{(n-m)*}(R\otimes S) - T_+^s\otimes T_-^s, 
d_+^{-2sm}(F^m)_*(\widehat\varphi \Theta)\right\rangle\\
& =: & \left\langle d_+^{-2s(n-2m)} F^{(n-2m)*}(T) - T_+^s\otimes T_-^s, 
\Phi\right\rangle
\end{eqnarray*}
where $T:= d_+^{-2sm} F^{m*}(R\otimes S)$ and $\Phi:=d_+^{-2sm} (F^m)_*(\widehat \varphi\Theta)$.

Hence, $T$ has support in a small neighbourhood ${\cal U}$ of the filled
Julia set $\K^F_+=\K_+\times \K_-$ of $F$
and
$\Phi$ is a smooth form with support in a small neighbourhood ${\cal V}$ of $\K_-^F=\K_-\times\K_+$. 
Moreover, since $m$ is big and $\varphi$ is p.s.h. on $U_2\cap V_2$, 
$\ddc\Phi\geq 0$ in a neighbourhood ${\cal U}'\Supset {\cal U}$ of $\K_+^F$. Putting 
$\widehat\omega(z,w):=\omega(z)$,
we have $-\widehat\omega\leq \ddc \widehat\varphi\leq\widehat\omega$ since $\|\varphi\|_{{\cal C}^2}=1$.
It follows that 
$$-d_+^{-2sm}(F^m)_*(\widehat\omega\wedge \Theta) \leq \ddc\Phi
\leq d_+^{-2sm}(F^m)_*(\widehat\omega\wedge \Theta).$$  

The choice of $\cal W$, $\cal U$, $\cal V$, ${\cal U}'$ and $m$ does not depend on $\varphi$ and $n$. 
Lemma 3.3 and  Proposition 2.1 applied to $F$, $T$ and $\Phi$ imply
$$\left\langle d_+^{-2s(n-2m)} F^{(n-2m)*}(T) - T_+^s\otimes T_-^s, 
\Phi\right \rangle \leq c'd_+^{-n}$$
and
$$\left|\left\langle d_+^{-2s(n-2m)} F^{(n-2m)*}(T) - T_+^s\otimes T_-^s, 
\Phi\right \rangle \right|\leq c'_Td_+^{-n}.$$
The desired inequalities of the proposition follow.
Since every smooth test 
function on $\P^k$ can be written as a difference of smooth functions p.s.h. on $U_2\cap V_2$,
these inequalities imply that 
$d_+^{-2sn}f^{n*}(R)\wedge (f^n)_*(S)\rightarrow\mu$.
\hfill $\square$

\begin{corollary} The Green measure of $F$ is equal to $\mu\otimes\mu$.
\end{corollary}
\begin{proof}
Let $R$ and $S$ be as in Proposition 3.1 such that $\supp(R\otimes S)\cap \{t=0\}\subset I^F_+$
and $\supp(S\otimes R)\cap \{t=0\}\subset I^F_-$. Proposition 3.1 implies that the Green measure of 
$F$ is equal to
\begin{eqnarray*}
\lefteqn{\lim d_+^{-4sn} F^{n*}(R\otimes S)\wedge (F^n)_*(S\otimes R)}\\
 & = & \lim d_+^{-4sn}\left[f^{n*}(R)\otimes (f^n)_*(S)\right] \wedge  \left[(f^n)_*(S)\otimes f^{n*}(R)\right]\\
& = & \lim \left[d_+^{-2sn}f^{n*}(R)\wedge (f^n)_*(S)\right] \otimes 
\left[d_+^{-2sn}f^{n*}(R)\wedge (f^n)_*(S)\right]\\
& = &\mu\otimes\mu.
\end{eqnarray*}
\end{proof}

\section{Speed of mixing}

In this section, we give the proof of Theorem 1.1. Fix a domain $D$ in $\C^k$ containing $\K:=\K_+\cap\K_-$. 
Observe that $\varphi$ and $\psi$ can be written as
differences of smooth functions strictly p.s.h. on a neighbourhood of $\overline D$.
Hence, we can assume that $\ddc\varphi\geq \omega$ on $D$, $\ddc\psi\geq\omega$ 
on $D$ and $\|\varphi\|_{{\cal C}^2}\leq M$, $\|\psi\|_{{\cal C}^2}\leq M$ for some fixed constant
$M>0$. The constants $c$, $A$, $c'$ below
do not depend on $\varphi$ and $\psi$. 

It is sufficient to prove Theorem 1.1 for $n$ even. So we have to prove that
\begin{eqnarray}
\left|\left\langle\mu, (\varphi\circ f^n)(\psi\circ f^{-n})\right\rangle 
- \langle \mu,\varphi\rangle \langle \mu, \psi 
\rangle \right|\leq c d_+^{-n}.
\end{eqnarray}
Observe that, 
since $\mu$ is invariant, the left hand side of (4) 
does not change if we add to $\varphi$ or to $\psi$ a constant.
Consequently, one only need to check for a constant $A$ that 
\begin{eqnarray}
\left\langle\mu, (\varphi\circ f^n +A)(\psi\circ f^{-n}+A)\right\rangle 
- \langle \mu,\varphi+A\rangle \langle \mu, \psi +A
\rangle \leq c d_+^{-n}
\end{eqnarray}
and
\begin{eqnarray}
\left\langle\mu, (\varphi\circ f^n -A)(-\psi\circ f^{-n}+A)\right\rangle 
- \langle \mu,\varphi-A\rangle \langle \mu, -\psi +A
\rangle \leq c d_+^{-n}.
\end{eqnarray}
We choose $A>0$ big enough so that $\phi(z,w):=(\varphi(z)+A)(\psi(w)+A)$ and 
$\phi'(z,w):=(\varphi(z)-A)(-\psi(w)+A)$ are p.s.h. on $D\times D$. 
This allows to apply Proposition 3.1 to the automorphism $F$ and to the test functions $\phi$, $\phi'$.
We will check (5). The estimate (6) can be proved in the same way.

Define $T_1:=T_+^s \otimes T_-^s$. 
Since $F^*(T_1)=d_+^{2s} T_1$ and $T_\pm$ have continuous 
potentials in $\C^k$, we get the following identities
\begin{eqnarray*}
\left\langle \mu, (\varphi\circ f^n+A)(\psi\circ f^{-n}+A) \right\rangle & = & 
\left\langle T_+^s\wedge T_-^s, (\varphi\circ f^n+A)
(\psi\circ f^{-n}+A)\right\rangle \\
& = & \left\langle T_1 \wedge [\Delta], \phi\circ F^n \right\rangle \\
& = & \left\langle d_+^{-4sn+2sm} F^{(2n-m)*}(T_1)\wedge [\Delta], \phi\circ F^n\right\rangle \\
& = & \left\langle d_+^{-4sn+2sm} F^{(n-m)*}(T_1)\wedge (F^n)_*[\Delta], \phi\right\rangle  \\
& =: & \left\langle d_+^{-4sn+4sm} F^{(n-m)*}(T_1)\wedge (F^{n-m})_*T_2, \phi\right\rangle
\end{eqnarray*}
where $T_2:=d^{-2sm}(F^m)_*[\Delta]$ and $m$ is a fixed, but sufficiently large integer.
By Lemma 3.2, $T_2$ has support 
in a small neighbouhood ${\cal V}$ of $\K_-^F$. 

Using a regularization of currents, we may find smooth currents $T'_1$ and $T_2'$ with support in small
neighbourhoods $\cal U$ of $\K^F_+$ and $\cal V$ of $\K^F_-$ respectively, so that
\begin{eqnarray*}
\lefteqn{\left \langle d_+^{-4sn+4sm} F^{(n-m)*}(T_1) \wedge (F^{n-m})_*T_2 ,\phi\right \rangle -}\\
&& 
-\left \langle d_+^{-4sn+4sm} F^{(n-m)*}(T_1') \wedge (F^{n-m})_*T_2', \phi\right \rangle \leq d_+^{-n} .
\end{eqnarray*}
The currents $T_1'$ and $T_2'$ depend on $n$. The choice of $m$, $\cal U$, $\cal V$ 
depends only on $D$ and $f$ with ${\cal U}\cap{\cal V}\Subset D\times D$.

Since $\langle\mu,\varphi+A\rangle\langle \mu,\psi+A\rangle=\langle \mu\otimes\mu,\phi\rangle$, we only have 
to check that
\begin{eqnarray*}
\left \langle d_+^{-4sn+4sm} F^{(n-m)*}(T_1') \wedge (F^{n-m})_*T_2'
-\mu\otimes\mu, \phi\right \rangle \leq c'd_+^{-n} .
\end{eqnarray*}
This inequality follows directly from Corollary 3.4 and Proposition 3.1 applied to $F$ and $\phi$. 
Hence, the proof of Theorem 1.1 is complete.
\\

\noindent
{\bf Acknowledgments.} This article is written during my visit to the Humboldt University of Berlin. I
would like to thank the Alexander 
von Humboldt foundation 
for its support and Professor J\"urgen Leiterer for his help and
his great hospitality.

Tien-Cuong Dinh,
Math\'ematique - B\^at. 425, UMR 8628, 
Universit\'e Paris-Sud, 91405 Orsay, France. 
E-mail: TienCuong.Dinh@math.u-psud.fr.
\end{document}